%&amstex          
\input amstex\documentstyle{amsppt}  
\pagewidth{12.5cm}\pageheight{19cm}\magnification\magstep1
\topmatter
\title{Precuspidal families and indexing of Weyl group representations}\endtitle
\author G. Lusztig\endauthor
\address{Department of Mathematics, M.I.T., Cambridge, MA 02139}\endaddress
\thanks{Supported by NSF grant DMS-2153741}\endthanks
\endtopmatter   
\document
\define\ovm{\overset\smile\to}
\define\occ{\ovm{\cc}}

\define\bx{\boxed}

\define\sneq{\subsetneqq}

\define\hW{\hat W}

\define\pe{\perp}
\define\si{\sim}

\define\sqc{\sqcup}

\define\bA{\bar A}

\define\bX{\bar X}
\define\bZ{\bar Z}
\define\lb{\linebreak}

\define\bin{\binom}
\define\op{\oplus}
   
\define\part{\partial}
\define\emp{\emptyset}

\define\m{\mapsto}
\define\do{\dots}

\define\sub{\subset}    
\define\bxt{\boxtimes}
\define\T{\times}
\define\ti{\tilde}
\define\nl{\newline}
\redefine\i{^{-1}}

\define\un{\underline}
\define\ucf{\un{\cf}}
\define\uocc{\un{\occ}}

\define\ot{\otimes}

\define\Hom{\text{\rm Hom}}

\define\Ind{\text{\rm Ind}}

\define\sgn{\text{\rm sgn}}

\define\he{\heartsuit}

\define\g{\gamma}
\redefine\d{\delta}
\define\e{\epsilon}
\define\et{\eta}
\define\io{\iota}

\define\p{\pi}

\define\r{\rho}
\define\s{\sigma}
\redefine\t{\tau}

\redefine\l{\lambda}

\define\x{\xi}

\redefine\G{\Gamma}
\redefine\D{\Delta}

\define\Si{\Sigma}
\define\Th{\Theta}

\define\Ph{\Phi}

\redefine\ss{\bold s}

\redefine\AA{\bold A}

\define\CC{\bold C}

\define\FF{\bold F}

\define\NN{\bold N}

\define\QQ{\bold Q}

\define\ZZ{\bold Z}

\define\cc{\Cal C}

\define\ce{\Cal E}
\define\cf{\Cal F}

\define\cl{\Cal L}

\define\car{\Cal R}

\define\ct{\Cal T}

\define\cz{\Cal Z}
\define\cx{\Cal X}

\define\fF{\frak F}

\define\tu{\ti u}

\define\tS{\ti S}

\head Introduction\endhead
\subhead 0.1\endsubhead
Let $W$ be the Weyl group of a connected reductive group $G$ over
$\CC$ and let $\hW$ be the set of isomorphism classes of irreducible
$\QQ[W]$-modules. In \cite{L79}, \cite{L82}, a partition of $\hW$ into
subsets called {\it families} was defined.
Let $\Ph(W)$ be the set of families of $\hW$. In \cite{L79}, \cite{L84}
we have attached to each $c\in\Ph(W)$ a certain finite group $\G_c$. 

Let $\car(W)$ be the $\CC$-vector space  with basis $\hW$.
For any $c\in\Ph(W)$ we denote by $\car_c$ the subspace
of $\car(W)$ spanned by the basis elements in $c$.
Let $E_c$ be the special representation in $c$ viewed as an element of
$\car_c$. On the other hand a set $Con_c$ of {\it constructible
representations} of $W$ wes associated to $c$ in \cite{L82}. Any
$\r\in Con_c$ can be viewed as an element of $\car_c$ (which is not
contained in the obvious basis of $\car_c$, if $\G_c\ne1$). In
\cite{L87} to each $\r\in Con_c$ we have attached a subgroup $H_\r$ of
$\G_c$. We can attach to $E_c$ the pair of subgroups $(\{1\},\G_c)$ of
$\G_c$. It is remarkable that one can interpolate between
$[(H_\r,H_\r),\r]$ ($\r\in Con_c$) on the one hand and $[(\{1\},\G_c),E_c]$
on the other hand. The resulting objects $[(\G',\G''),e_{\G',\G''}]$
are such that $\G',\G''$ are subgroups of $\G_c$ with $\G'$ normal in
$\G''$ and the elements $e_{\G',\G''}$ form a basis $b_c$ of $\car_c$
which is related to the standard basis of $\car_c$ by an upper
triangular matrix with entries in $\NN$ and with $1$ on diagonal; thus
$b_c$ is in canonical bijection with $c$ and also in bijection with
the set $X(\G_c)$ of pairs $(\G',\G'')$ as above.
Thus the representations in $c$ can be indexed by the pairs of groups
in $X_{\G_c}$.
This has been done in \cite{L19},\cite{L20},\cite{L22},\cite{L23}, but
not in a uniform manner.

In this paper we do the same in a uniform manner.
We define $X_{\G_c}$ by an inductive procedure which
starts with a subset $x_{\G_c}$ of $X_{\G_c}$ which is much smaller
than $X_{\G_c}$. To describe $x_{\G_c}$ it suffices to consider the
case where $W$ is irreducible and $c$ is a cuspidal family (see
4.4) of $W$. 
The pairs in $x_{\G_c}$ are then essentially
described in terms of the various ways in which $c$ can be obtained
by $J$-induction (see 4.2) from certain families (which we call
{\it precuspidal}, see 4.5, 4.6) in various proper parabolic subgroups
of $W$.

In \S5 an equivalent description of $x_{\G_c}$ in terms of the various
ways in which the special unipotent class defined by $c$ can be obtained
by the induction procedure of \cite{LS79} from unipotent classes in
Levi subgroups of proper parabolic subgroups of $G$ is given.

\subhead 0.2\endsubhead
{\it Notation.} The number of elements in a finite set $X$ is denoted by
$|X|$. Let $\FF$ be the field $\ZZ/2\ZZ$.
An {\it interval} in $\NN$ is a subset of $\NN$ of the form \lb
$[a,b]=\{c\in\NN;a\le c\le b\}$ where $a\le b$ are in $\NN$;
we write $a<<b$ whenever $b-a\ge2$.
For $\d\in[0,1]$ and $H\sub\ZZ$ we set $H^\d=H\cap(\d+2\ZZ)$.
For an element $u$ in a group $G$ we denote by $Z_G(u)$
the centralizer of $u$ in $G$; if $G$ is an algebraic group
let $A_G(u)$ be the group of components of $Z_G(u)$.

\head 1. Subspaces of an $\FF$-vector space\endhead
\subhead 1.1\endsubhead
Let $V$ be the  $\FF$-vector space with basis
$\{e_i;i\in\NN_{>0}\}$. For $a\le b$ in $\NN_{>0}$ we set
$e_{[a,b]}=e_a+e_{a+1}+\do+e_b\in V$.

Let $\cf(V)$ be the set of subspaces $E$ of $V$
such that $E$ has a basis \lb $\{e_{[a_k,b_k]};k=1,\do,r\}$
where $[a_1,b_1],[a_2,b_2],\do,[a_r,b_r]$
are intervals in $\NN_{>0}$ for which (i)-(iii) below hold.

(i) For any $k\in[1,r]$ we have $a_k=b_k\mod2$.

(ii) If $k\in[1,r]$ and $a_k<c<b_k$ with $c-a_k=1\mod2$ then there
exists $k'\in[1,r]$ such that $a_k<a_{k'}\le c\le b_{k'}<b_k$.

(iii) For any $k\ne k'$ in $[1,r]$ we have $b_k<<a_{k'}$,
or $b_{k'}<<a_k$, or $a_k<a_{k'}\le b_{k'}<b_k$, or
$a_{k'}<a_k\le b_k<b_{k'}$.
\nl
Note that for $E\in\cf(V)$, the intervals $[a_k,b_k]$ as
above are uniquely determined by $E$, so that for any
$j\in\NN_{>0}$ we can define $f_j(E)\in\NN$ to be the
number of $k\in[1,r]$ such that $j\in[a_k,b_k]$. We set
$$\e(E)=\sum_{j\in\NN_{>0}}(1/2)f_j(E)(f_j(E)+1)e_j\in V.$$
We obtain a function $\e:\cf(V)@>>>V$.

\subhead 1.2\endsubhead
Let $D\in\NN$. Let $V_D$ be the subspace of $V$ with basis
$\{e_i;i\in[1,D]\}$. Assuming that $D\ge2$ and $j\in[1,D]$,
let $U_{D,j}$ be the subspace of $V_D$ with basis
consisting of $e_j$ and of

$(*)$:

$e_3,e_4,\do,e_D$, if $j=1$;

$e_1,\do,e_{j-2},e_{[j-1,j+1]},e_{j+2},\do,e_D$, if $1<j<D$;

$e_1,e_2,\do,e_{D-2}$, if $j=D$.
\nl
Let $C_j:U_{D,j}@>>>V_{D-2}$ be the surjective linear map 
which carries $e_j$ to $0$ and carries the vectors in $(*)$
(in the order written) to the vectors $e_1,e_2,\do,e_{D-2}$
(in the order written).

Following \cite{L22}, \cite{L23}
(or \cite{L19} in the case where $D$ is even)
we define a collection $\cf(V_D)$ of subspaces of $V_D$ by
induction on $D$.
If $D=0$, $\cf(V_D)$ consists of $\{0\}$. If $D=1$, $\cf(V_D)$
consists of $\{0\}$ and of $V_D$. Assume now that $D\ge2$. We
say that a subspace $E$ of $V_D$ is in $\cf(V_D)$ if either
$E=0$ or if there exist $j\in[1,D]$ and $E'\in\cf(V_{D-2})$
such that $E=C_j\i(E')$.

The following result is stated in \cite{L22} and proved in
the case where $D$ is even in \cite{L19} and in the case where $D$
is odd in \cite{L23}.

(a) $\cf(V_D)=\{E\in\cf(V);E\sub V_D\}$.

\subhead 1.3\endsubhead
We define a map $u:V@>>>\ZZ$ as follows. If $x\in V$ we can
write $x$ uniquely in the form
$$x=e_{[a_1,b_1]}+e_{[a_2,b_2]}+\do+e_{[a_r,b_r]}\tag a$$
where $1\le a_1\le b_1<<a_2\le b_2<<\do<<a_r\le b_r$. We set
$$u(x)=\sum_{s\in[1,r];a_s+b_s=1\mod2}(-1)^{a_s}\in\ZZ.$$
This defines $u$. Let ${}^0V=u\i(0)$.

For any $D\ge0$ we set ${}^0V_D={}^0V\cap V_D$.
From \cite{L19, 1.16} and its proof,

(b) {\it for any even $D\ge0$, the map $\e:\cf(V)@>>>V$
restricts to a bijection $\e_D:\cf(V_D)@>\si>>{}^0V_D$
such that $\e_D(E)\in E$ for any $E\in\cf(V_D)$.}

Since $\cf(V)=\cup_{D\ge0,even}\cf(V_D)$ and
${}^0V=\cup_{D\ge0,even}{}^0V_D$, it follows that

(c) {\it $\e$ defines a bijection $\cf(V)@>\si>>{}^0V$; moreover,
we have $\e(E)\in E$ for any $E\in\cf(V)$.}

Arguments similar to those in \cite{L19, 1.16} show that

(d) {\it statement (b) remains valid if ``even'' is replaced
by ``odd''.}
\nl
Let $Z$ be the $\FF$-vector space with basis
$\{g_i;i\in\NN\}$. For $z\in Z$ let
${}^0\tu(z)$ (resp. ${}^1\tu(z)$) be the number of
even (resp. odd) $i\in\NN$ such that $g_i$
appears with nonzero coefficient in $z$. We define
$\tu:Z@>>>\ZZ$ by $\tu(z)={}^0\tu(z)-{}^1\tu(z)$.
For $1\le a\le b$ we set $g_{[a,b]}=g_a+g_{a+1}+\do+g_b$.

Let $\bZ$ be the subspace of $Z$ consisting of the elements
$\sum_ia_ig_i$ with $a_i\in\FF$, $\sum_ia_i=0$. We define an
isomorphism $\x:V@>>>\bZ$ by $\x(e_i)=g_{i-1}+g_i$ for all
$i\in\NN_{>0}$. We show:

(e) If $x\in V$, then $u(x)=-\tu(\x(x))$.
\nl
We write $x$ as in (a). We have 
$$\x(x)=g_{[a_1-1,b_1]}+g_{[a_2-1,b_2]}+\do+g_{[a_r-1,b_r]}.$$
Since $0\le a_1-1<b_1\le a_2-1<b_2\le\do\le a_r-1<b_r$, we have
(for $\d=0,1$):
$$\align&{}^\d\tu(\x(x))=\sum_k{}^\d\tu(g_{[a_k-1,b_k]}=
\sum_{k;a_k-b_k\in\ZZ^0}(b_k-a_k+2)/2 \\&+
\sum_{k;a_k\in\ZZ^\d,b_k\in\ZZ^{1-\d}}(b_k-a_k+1)/2+
\sum_{k;a_k\in\ZZ^{1-\d},b_k\in\ZZ^\d}(b_k-a_k+3)/2\endalign$$

hence    
$$\align&\tu(\x(x))=\sum_{k;a_k\in\ZZ^0,b_k\in\ZZ^1}
(b_k-a_k+1)/2+
\sum_{k;a_k\in\in\ZZ^1,b_k\in\ZZ^0}(b_k-a_k+3)/2\\&-
\sum_{k;a_k\in\ZZ^0,b_k\in\ZZ^1}(b_k-a_k+3)/2-
\sum_{k;a_k\in\ZZ^1,b_k\in\ZZ^0}(b_k-a_k+1)/2\\&
=\sum_{k;a_k\in\ZZ^0,b_k\in\ZZ^1}(-1)+
\sum_{k;a_k\in\ZZ^1,b_k\in\ZZ^0}1=-u(x).\endalign$$

This proves (e).

Now let $D\ge1$. Let $Z_D$ be the subspace of $Z$ spanned by
$\{g_i;i\in[0,D]\}$ and let $\bZ_D=Z_D\cap\bZ$. Now $\x$
restricts to a bijection $V_D@>>>\bZ_D$ and (by (e)) this
restricts to a bijection

(f) ${}^0V_D@>\si>>{}^0\bZ_D$
\nl
where ${}^0\bZ_D=\{y\in\bZ_D;\tu(y)=0\}$.
Let $Z'_D$ be the set of all $H\sub[0,D]$ such that $|H^0|=|H^1|$.
Clearly, $H\m\sum_{i\in H}g_i$ is a bijection
$Z'_D@>si>>{}^0\bZ_D$. Let $Z''_D$ be the set of $H\sub[0,D]$
such that $|H|=(D+2)/2$ (if $D$ is even) or
$|H|=(D+1)/2$ (if $D$ is odd).
We have a bijection $Z'_D@>\si>>Z''_D$ given by
$H\m([0,D]^0-H^0)\cup H^1$.
Using this and (f) we see that $|{}^0V_D|=|Z''_D|$. Thus:

{\it if $D$ is even then
$|{}^0V_D|=\bin{D+1}{(D+2)/2}=\bin{D+1}{D/2}$};

{\it if $D$ is odd then $|{}^0V_D|=\bin{D+1}{(D+1)/2}$.}

Now if $D$ is odd, $Z'_D$ has a fixed point free involution
$H\m [0,D]-H$.
\nl
This corresponds to the fixed point free involution
$z\m z+\sum_{i\in[0,D]}g_i$ that is, $z\m z+\x(e_1+e_3+\do+e_D)$
of ${}^0\bZ_D$. 
This corresponds under $\x$ to the fixed point free involution

(g) $\Th:x\m x+(e_1+e_3+\do+e_D)$ 
\nl
of ${}^0V_D$; in particular:

(h) ${}^0V_D$ is stable under $\Th$.

\subhead 1.4\endsubhead
We define a symplectic form $(,):V\T V@>>>\FF$ by
$(e_i,e_j)=1$ if $i-j=\pm1$ and $(e_i,e_j)=0$ if $i-j\ne\pm1$.
Let $D\ge0$. From the inductive definition of $\cf(V_D)$
(see 1.2) we see that 

(a) {\it any $E\in\cf(V_D)$ satisfies $(E,E)=0$.}

\subhead 1.5\endsubhead
We have $V=V^0\op V^1$ where $V^0$ (resp. $V^1$)
is the subspace spanned by $e_2,e_4,e_6,\do$
(resp. by $e_1,e_3,e_4,\do$).
For any $D\ge0$ we have $V_D=V_D^0\op V_D^1$ where
$V_D^0=V_D\cap V^0,V_D^1=V_D\cap V^1$.

Assuming that $D\ge2$ and $j\in[1,D]$, we have
$U_{D,j}=U_{D,j}^0\op U_{D,j}^1$
where $U_{D,j}^0=U_{D,j}\cap V_D^0$,
$U_{D,j}^1=U_{D,j}\cap V_D^1$.

If $j$ is odd, we have $U_{D,j}^1=V_D^1$.
If $j$ is even, (setting $D^-=D$ if $D$ is odd,
$D^-=D-1$ if $D$ is even), $U_{D,j}^1$ is the subspace of
$V_D^1$ with basis consisting of

$(*)$:

$e_1,e_3,\do,e_{j-3},e_{j-1}+e_{j+1},e_{j+3},e_{j+5},\do,
e_{D^-}$ (if $1<j<D$);

$e_1,e_3,\do,e_{D^--2}$, if $j=D$.

Now $C_j:U_{D,j}@>>>V_{D-2}$ induces a surjective linear map
$C_j^1:U_{D,j}^1@>>>V_{D-2}^1$.

When $j$ is odd, $C^1_j:U_{D,j}^1=V_D^1@>>>V_{D-2}^1$
carries $e_j$ to $0$ and carries the vectors
$e_1,e_3,\do,e_{j-2},e_{j+2},e_{j+4},\do, e_{D^-}$
(in the order written) to the vectors
$e_1,e_3,\do,e_{D^--2}$ (in the order written).

When $j$ is even, $C^1_j:U_{D,j}^1@>>>V_{D-2}^1$ 
carries the vectors in $(*)$ (in the order written) to the vectors 
$e_1,e_3,\do,e_{D^--2}$ (in the order written).

We define a collection $\occ(V_D^1)$ of pairs
$(\cl\sub\cl')$ of subspaces of $V_D^1$ by induction on $D$.
If $D=0$, $\occ(V_D^1)$ consists of $(\{0\}\sub\{0\})$.
If $D=1$, $\occ(V_D^1)$
consists of $(V_D^1\sub V_D^1)$ and $(\{0\}\sub V_D^1)$.
Assume now that $D\ge2$. We say that a pair $(\cl\sub\cl')$ of
subspaces of $V_D^1$ is in $\occ(V_D^1)$ if either
$(\cl\sub\cl')=(0\sub V_D^1)$ or if there exist
$j\in[1,D]$ and $(\cl_1\sub\cl'_1)\in\occ(V_{D-2}^1)$
such that $\cl=(C_j^1)\i(\cl_1)$, $\cl'=(C_j^1)\i(\cl'_1)$.

Using induction on $D$ we see that

(a) {\it if $D$ is odd and $(\cl\sub\cl')\in\occ(V_D^1)$ then
$e_1+e_3+\do+e_D\in\cl'$.}

\subhead 1.6\endsubhead
If $E\in\cf(V_D)$ we have $E=E^0\op E^1$ where
$E^0=E\cap V^0_D$, $E^1=E\cap V^1_D$. 

To any $E\in\cf(V_D)$ we associate the subspace
$(E^0)^!=\{x\in V_D^1;(x,E^0)=0\}$ of $V_D^1$.
Using 1.4(a) we see that

(a) $E^1\sub(E^0)^!$.
\nl
The following result is proved by induction on $D$.

(b) {\it The map $E\m(E^1\sub(E^0)^!)$ is a well defined bijection
$\Pi_D:\cf(V_D)@>\si>>\occ(V_D^1)$.}
\nl
(From this one can deduce an alternative proof of 1.5(a).)

In \cite{L22, 1.9}, a set also denoted by $\occ(V_D^1)$ is
defined in a way different from the way it is defined here;
however the two definitions agree, as a consequence of (b)
and \cite{L22, 1.9(b)}, \cite{L23}.

\subhead 1.7\endsubhead
In this subsection we assume that $D\ge1$ is odd.
We set $\et_D=e_1+e_3+\do+e_D\in V_D$.
Let $V'_D=V_D/\FF\et_D$. Let $\p:V_D@>>>V'_D$ be the obvious
projection. We have $V'_D=V'_D{}^0\op V'_D{}^1$ where
$V'_D{}^0=\p(V_D^0),V'_D{}^1=\p(V_D^1)$.
Now $(,)$ induces a nondegenerate symplectic form
$(,)':V'_D\T V'_D@>>>\FF$.

We define $\l:\cf(V_{D-1})@>>>\{\text{set of subspaces of} V'_D\}$
by $E\m\p(E)$. Since $V_{D-1}\cap\ker\p=0$, $\l$ is injective.
We denote by $\ucf(V'_D)$ its image. Thus $\l$ can be viewed as a
bijection $\cf(V_{D-1})@>\si>>\ucf(V'_D)$.
If $\ce\in\ucf(V'_D)$ we have $\ce=\ce^0\op\ce^1$ where
$\ce^0=\ce\cap V'_D{}^0$, $\ce^1=\ce\cap V'_D{}^1$; we set
$(\ce^0)^!=\{x'\in V'_D{}^1;(x',\ce^0)'=0\}$. Note that
$\ce^1\sub(\ce^0)^!$. We define
$$\l':\ucf(V'_D)@>>>\{\text{set of pairs of subspaces of} V'_D{}^1\}$$
by $\ce\m(\ce^1\sub(\ce^0)^!)$. We denote by $\uocc(V'_D{}^1)$ the
image of $\l'$. Thus $\l'$ can be viewed as a
surjective map $\ucf(V'_D)@>\si>>\uocc(V'_D{}^1)$.
This map is in fact a bijection: assume that $\ce,\ce'$ in
$\ucf(V'_D)$ satisfy $\ce^1=\ce'{}^1$, $(\ce^0)^!=(\ce'{}^0)^!$;
using the nondegeneracy of $(,)'$ we deduce $\ce^0=\ce'{}^0$
hence $\ce=\ce'$, as desired.

Let ${}^0V'_D=\p({}^0V_D)$. 

We define $\e':\ucf(V'_D)@>>>{}^0V'_D$ by $\e'(\p(E))=\p(\e_{D-1}(E))$
with $E\in\cf(V_{D-1})$.
If $E,E'$ in $\cf(V_{D-1})$ satisfy $\e'(\p(E))=\e'(\p(E'))$ the
$\e_{D-1}(E),\e_{D-1}(E')$ are in the same fibre of
$\p:{}^0V_D@>>>{}^0V'_D$ hence are either equal or their difference
is $\et_D$; the second possibility cannot occur since
$\e_{D-1}(E),\e_{D-1}(E')$ are both in $V_{D-1}$. We see that
$\e_{D-1}(E)=\e_{D-1}(E')$ so that $E=E'$, see 1.3(b).
Thus $\e'$ is injective.

From 1.3(h) we see that
$|{}^0V'_D|=(1/2)|{}^0V_D|=(1/2)\bin{D+1}{(D+1)/2}$.
From 1.3 we have
$|\ucf(V'_D)|=|\cf(V_{D-1})|=|{}^0V_{D-1}|=\bin{D}{(D-1)/2}$.
Since $\bin{D}{(D-1)/2}=(1/2)\bin{D+1}{(D+1)/2}$ we see that
$|\ucf(V'_D)|=|{}^0V'_D|$. Since $\e'$ is an injective map between
finite sets with the same cardinal, we see that:

(a) {\it $\e'$ is a bijection.}

\subhead 1.8\endsubhead
In this subsection we assume that $D\ge1$ is odd.
Assuming that $D\ge3$ and $j\in[1,D]$ we set
$U'_{D,j}=\p(U_{D,j})$. We have
$U'_{D,j}=U'_{D,j}{}^0\op U'_{D,j}{}^1$ where
$U'_{D,j}{}^0=\p(U_{D,j}^0)$, $U'_{D,j}{}^1=\p(U_{D,j}^1)$.
Since $C_j(\et_D)=\et_{D-2}$, the linear maps
$C_j:U_{D,j}@>>>V_{D-2}$, $C_j^1:U_{D,j}^1@>>>V_{D-2}^1$
induce linear maps $C'_j:U'_{D,j}@>>>V'_{D-2}$,
$C'_j{}^1:U'_{D,j}{}^1@>>>V'_{D-2}{}^1$.

We define a collection $\cf(V'_D)$ of subspaces
of $V'_D$ by induction on $D$. If $D=1$, $\cf(V'_D)$
consists of $\{0\}$. Assume now that $D\ge3$. We say that
a subspace $\ce$ of $V'_D$ is in $\cf(V'_D)$ if either $\ce=0$
or if there exist $j\in[1,D-1]$ and $\ce'\in\cf(V'_{D-2})$
such that $\ce=(C'_j)\i(\ce')$. (Note that in this definition
$j$ is not allowed to be $D$.)

We define a collection $\occ(V'_D{}^1)$ of pairs
$(\cl\sub\cl')$ of subspaces of $V'_D{}^1$ by induction on $D$.
If $D=1$, $\occ(V'_D{}^1)$ consists of $(\{0\}\sub V'_D{}^1)$.
Assume now that $D\ge3$. We say that a pair $(\cl\sub\cl')$ of
subspaces of $V'_D{}^1$ is in $\occ(V'_D{}^1)$ if either
$(\cl\sub\cl')=(0\sub V'_D{}^1)$ or if there exist
$j\in[1,D-1]$ and $(\cl_1\sub\cl'_1)\in\occ(V'_{D-2}{}^1)$
such that $\cl=((C'_j)^1)\i(\cl_1)$, $\cl'=((C'_j)^1)\i(\cl'_1)$.
(Again, $j$ is not allowed to be $D$.)

From the definitions we see that

$\ucf(V'_D)=\cf(V'_D)$.
\nl
We use that
under the isomorphism $V_{D-1}@>\si>>V'_D$ 
induced by $e_i\m e_i$, the operators $C_j$ ($j\in[1,D-1]$) used
to define the left hand side correspond to the 
operators $C'_j$ ($j\in[1,D-1]$) used to define the right hand side.
(Note that $C'_j$ with $j=D$ is not used and $C_j$
with $j=D$ is 
not defined.)

Similarly we have $\uocc(V'_D{}^1)=\occ(V'_D{}^1)$.
Hence, using 1.7(a),

(a) {\it $\e'$ in 1.7 can be regarded as a bijection $\cf(V'_D)@>>>{}^0V'_D$.}

and

(b) {\it $\l'$ in 1.7 can be regarded as a bijection
$\cf(V'_D)@>\si>>\occ(V'_D{}^1)$.}

\head 2. The sets $x_\G,X_\G$\endhead
\subhead 2.1\endsubhead
For any $n\ge1$ we denote by $S_n$ the symmetric group
consisting of all permutations of $[1,n]$.
If $n\ge2$ we identify $S_{n-1}$ with the subgroup of $S_n$
consisting of permutations of $[1,n]$ which keep $n$ fixed.
Thus we have $S_1\sub S_2\sub S_3\sub\do$.
We denote by $\D_8$ the centralizer in $S_4$ of the permutation
$1\m2,2\m1,3\m4,4\m3$ (a dihedral group; this is then also a
subgroup of $S_5,S_6$,etc.) We denote by $S_2S_2$ the subgroup of $S_4$
generated by the
transposition $1\m2,2\m1$ and by the transposition $3\m4,4\m3$
(this is a subgroup of $\D_8$). We denote by $\tS_2$ the subgroup of
$S_5$ generated by the transposition $4\m5,5\m4$. We denote by
$S_3S_2$ the centralizer in $S_5$ of that transposition.

\subhead 2.2\endsubhead
For any finite group $\G$ we denote by $\cz_\G$ the set of
pairs $(\G'\sub G'')$ of subgroups of $\G$ with $\G'$ normal
in $\G''$.

Let $\AA$ be the collection of finite groups consisting of
the following groups:

$V_D^1$ with  $D\in\NN$ even;  see 1.5.

$V'_D{}^1$ with $D\in\NN$ odd; see 1.7.

$S_n$ with $n\in[1,5]$;

$S'_2,S'_3$.
\nl
Here $S'_2$ (resp. $S'_3$) is another copy of $S_2$ (resp.
$S_3$). We view $S'_2,S_2$ as distinct objects of $\AA$.
We view $S'_3,S_3$ as distinct objects of $\AA$.
We view $V_0^1,V'_1{}^1,S_1$ as the same objects of $\AA$.
We view $V_2^1,V'_3{}^1,S_2$ as the same objects of $\AA$.
For $D$ odd, $D\ge5$, we view $V'_D{}^1$, $V_{D-1}^1$ as distinct
objects of $\AA$, although they are isomorphic as abstract groups.

To any $\G\in\AA$ with $|\G|>1$  we shall associate a subset
$x_\G$ of $\cz_\G$. In each case for each $(\G',\G'')$ in $x_\G$
we will describe

(a) the corresponding quotient $\G''/\G'$ (which turns out to be
again an object of $\AA$).

For $\G=V_D^1$ with $D\ge2$ even, $x_\G$ consists
of the pairs $(\G'_j,\G''_j)$ with $j\in[1,D]$ where

$\G''_j=U_{D,j}^1$ (see 1.5) and $\G'_j=0$ (if $j$ is even),
$\G'_j=\FF e_j$ (if $j$ is odd).
\nl
Note that $C_j^1$ (see 1.5) defines an isomorphism
$\G''_j/\G'_j@>>>V_{D-2}^1$.

For $\G=V'_D{}^1$ with $D\ge3$ odd, $x_\G$ consists of the
pairs $(\G'_j,\G''_j)$ with $j\in[1,D-1]$ where

$\G''_j=U'_{D,j}{}^1$ (see 1.8) and $\G'_j=0$ (if $j$ is even),
$\G'_j=\FF e'_j$ (if $j$ is odd).

Note that $C'_j{}^1$ (see 1.8) defines an isomorphism
$\G''_j/\G'_j@>>>V'_{D-2}{}^1$.

For $\G=S_2$, $x_\G$ consists of the pairs
$(S_1\sub S_1),(S_2\sub S_2)$; the corresponding quotients are $S_1$ and $S_1$.

For $\G=S_3$, $x_\G$ consists of the pairs $(S_1\sub S_2),
(S_3\sub S_3)$;
the corresponding quotients are $S_2$ and $S_1$.

For $\G=S'_2$, $x_\G$ consists of the pair $(S_2\sub S_2)$.
the corresponding quotient is  $S_1$.

For $\G=S'_3$, $x_\G$ consists of the pairs $(S_2\sub S_2),(S_3\sub S_3)$;
the corresponding quotients are $S_1$ and $S_1$.

For $\G=S_4$, $x_\G$ consists of the pairs
$(S_2S_2\sub\D_8),(S_2\sub S_2S_2),(S_3\sub S_3),(S_4\sub S_4)$;
the corresponding quotients are $S_2, S_2,S_1,S_1$.

For $\G=S_5$, $x_\G$ consists of the pairs
$(\tS_2\sub S_3S_2),(S_3\sub S_3S_2),(S_2S_2\sub\D_8),
(S_4\sub S_4),(S_5\sub S_5)$;
the corresponding quotients are $S_3,S_2,S_2,S_1,S_1$.

Note that if $\G$ is as above, we have $|\G'|=1$ for some $(\G',\G'')\in x_\G$ except when

$\G=S'_2$ (we then say that $\G$ is {\it anomalous},  or

$\G$ is one of $S'_3,S_4,S_5$ (we then say that $\G$ is
{\it terminal}).

\subhead 2.3\endsubhead
To any $\G\in\AA$ we shall associate a subset
$X_\G$ of $\cz_\G$ by induction on $|\G|$. If
$|\G|=1$, $X_\G$ consists of the pair $(\G\sub\G)$.
Assume now that $|\G|\ge2$. 

For any $(\G'\sub\G'')$ in $x_\G$ 
and any $(\G'_1\sub\G''_1)$ in $X_{\G''/\G'}$  
we define $\ti\G'_1,\ti\G''_1$ to be the inverse images of
$\G'_1,\G''_1$ under the quotient map $\G''@>>>\G''/\G'$.
(We have used 2.2(a).)
The pairs $(\ti\G'_1\sub\ti\G''_1)$ thus associated to various
$(\G'\sub\G'')$ in $x_\G$ 
and any $(\G'_1\sub\G''_1)$ in $X_{\G''/\G'}$ 
form a subset $(X_\G)_0$ of $\cz_\G$.
Note that for $(\G'_1\sub\G''_1)$ in $X_\G$  we have
$|\G''_1/\G'_1|<|\G|$; in particular, $(\G'_1\sub\G''_1)\ne(S_1\sub\G)$.

Let $Q(\G)$ be the set of subgroups $\G_1$ of $\G$ such that $(\G_1,\G_1)\in(X_\G)_0$.
Let $Q_*(\G)$ be the set of subgroups $\G_1\in Q(\G)$ such
that $\G_1$ is isomorphic to a product of groups in $\AA$;
the last condition is automatically satisfied
except when $\G$ is $S_4$ or $S_5$ in which case that
condition excludes $\G_1=\D_8$.

If $\G$ is not terminal we set $X_\G=(X_\G)_0\sqc\{(S_1\sub\G)\}$.

If $\G$ is terminal we set 
$X_\G=(X_\G)_0\sqc\sqc_{\G_1\in Q_*(\G)}\{(S_1\sub\G_1)\}$.

This completes the inductive definition of $X_\G$.
Note that for any
$(\G'_1\sub\G''_1)$ in $X_\G$,
$\G''_1/\G'_1$ is isomorphic to a product of groups in $\AA$.

From the definitions we see that:

If $\G=V_D^1$ with $D$ even, we have $X_\G=\occ(V_D^1)$
(see 1.5).

If $\G=V'_D{}^1$ with $D$ odd, we have
$X_\G=\occ(V'_D{}^1)$ (see 1.8).

If $\G=S_1$, $X_\G$ consists of $(S_1\sub S_1)$.

If $\G=S_2$, $X_\G$ consists of $(S_2\sub S_2),(S_1\sub S_2),(S_1\sub S_1)$.

If $\G=S_3$, $X_\G$ consists of $(S_3\sub S_3),(S_1\sub S_3),(S_2\sub S_2),(S_1\sub S_2),(S_1\sub S_1)$.

If $\G=S'_2$, $X_\G$ consists of $(S_2\sub S_2),(S_1\sub S_2)$.

If $\G=S'_3$, $X_\G$ consists of $(S_3\sub S_3),(S_1\sub S_3),(S_2\sub S_2),(S_1\sub S_2)$.

If $\G=S_4$, $X_\G$ consists of
$(S_4\sub S_4),(S_1\sub S_4),(\D_8\sub\D_8),
(S_2S_2\sub\D_8)$,
$(S_2S_2\sub S_2S_2),(S_2\sub S_2S_2),(S_1\sub S_2S_2)$,
$(S_3\sub S_3),(S_1\sub S_3),(S_2\sub S_2),(S_1\sub S_2)$.

If $\G=S_5$, $X_\G$ consists of
$(S_5\sub S_5),(S_1\sub S_5),(S_3S_2,S_3S_2)$,
$(S_3\sub S_3S_2),(\tS_2\sub S_3S_2),(S_1\sub S_3S_2)$,
$(S_4\sub S_4),(S_1\sub S_4),(\D_8\sub\D_8)$,
$(S_2S_2\sub\D_8)$,
$(S_2S_2\sub S_2S_2),(S_2\sub S_2S_2),(S_1\sub S_2S_2)$,
$(S_3\sub S_3),(S_1\sub S_3),(S_2\sub S_2),(S_1\sub S_2)$.

\subhead 2.4\endsubhead
For any finite group $\G$, we denote by $M(\G)$ the set of
$\G$-conjugacy of pairs $(x,\s)$ where $x\in\G$ and $\s$ is an
irreducible representation over $\CC$ of the centralizer
$Z_\G(x)$ of $x$ in $\G$.
Let $\CC[M(\G)]$ be the $\CC$-vector space with basis $M(\G)$. 
For $(\G'\sub\G'')\in\cz_\G$ let $\ss_{\G',\G''}:\CC[M(\G''/\G')]@>>>\CC[M(\G)]$
be the $\CC$-linear map defined in \cite{L20, 3.1}.
Now let $\G\in\AA$. For
$(\G'\sub\G'')\in X_\G$ we set $\r_{(\G'\sub\G'')}=\ss_{\G',\G''}(1,1)$
where $(1,1)$ is the element $(x,\s)\in M(\G''/\G')$ in which $x$ is the unit
element and $\s$ is the unit representation of $\G''/\G'$.
For example, $\r_{(S_1\sub\G)}=(1,1)\in M(\G)$.
Let  $M(\G)_0$ be the set of all $(x,\s)\in M(\G)$ such that $(x,\s)$ appears with
nonzero coefficient in $\r_{(\G',\G'')}$ for some $(\G'\sub\G'')\in X_\G$
and let $\CC[M(\G)_0]$ be the subspace of $\CC[M(\G)]$ spanned
by $M(\G)_0$.

The following result is a reformulation of results in \cite{L19}, \cite{L22}, \cite{L23}.
\proclaim{Theorem 2.5} (i) $\{\r_{(\G',\G'')};(\G'\sub\G'')\in X_\G\}$
is a $\CC$-basis of $\CC[M(\G)_0]$.

(ii) There is a unique bijection $j:M(\G)_0@>\si>> X_\G$ such that
for any $(x,\s)\in M(\G)_0$, $(x,\s)$ appears with coefficent $1$ in
$\r_{j(x,\s)}$.
\endproclaim
In the case where $\G$ is $V_D^1$ with $D$ even, we have
$M(\G)=V_D^1\op\Hom(V_D^1,\CC^*)=V_D^1\op\Hom(V_D^1,\FF)=V_D^1\op V_D^0$.
(The last equality comes by using $(,)$.) Hence $M(\G)=V_D$.
We have also $M(\G)_0={}^0V_D$. Then the bijection $j$ becomes a bijection
${}^0V_D@>>>\occ(V_D^1)$ or using the identification
$\occ(V_D^1)=\cf(V_D)$, a bijection ${}^0V_D@>>>\cf(V_D)$. This
coincides with the inverse of $\e_D$ in 1.3(b).

In the case where $\G$ is $V'_D{}^1$ with $D$ odd, we have
$M(\G)=V'_D{}^1\op\Hom(V'_D{}^1,\CC^*)=V'_D{}^1\op\Hom(V'_D{}^1,\FF)
=V'_D{}^1\op V'_D{}^0$
(The last equality comes by using $(,)')$. Hence $M(\G)=V'_D$.
We have also $M(\G)_0={}^0V'_D$. Then the bijection $j$ becomes a bijection
${}^0V'_D@>>>\occ(V'_D{}^1)$ or using the identification
$\occ(V'_D{}^1)=\cf(V'_D)$, a bijection ${}^0V'_D@>>>\cf(V'_D)$. This
coincides with the inverse of $\e'$ in 1.7(a).

\subhead 2.6\endsubhead
Let $\G\in\AA$.
Let $\le$ be the transitive relation on $M(\G)_0$ generated by
the relation for which $(x,\s),(x',\s')$ are related if
$(x,\s)$ appears with nonzero coefficient in $\r_{j(x',\s')}$.
The following result is a reformulation of results in \cite{L19}, \cite{L22}, \cite{L23}.
\proclaim{Theorem 2.7} $\le$ is a partial order on $M(\G)_0$. 
\endproclaim

\head 3. The sets $\bar x_\G,\bX_\G$\endhead
\subhead 3.1\endsubhead
To any $\G\in\AA$ with $|\G|>1$  we shall associate a subset
$\bar x_\G$ of $\cz_\G$.

We have $\bar x_\G=x_\G$ in all cases except when 
$\G=V'_D{}^1$ with $D\ge3$ odd, in which case $\bar x_G$
is defined like $x_\G$ but the condition $j\in[1,D-1]$ is
replaced by the condition $j\in[1,D]$.

\subhead 3.2\endsubhead
To any $\G\in\AA$ we shall associate a subset
$\bX_\G$ of $\cz_\G$ by induction on $|\G|$. If
$|\G|=1$, $\bX_\G$ consists of the pair $(\G\sub\G)$.
Assume now that $|\G|\ge2$. 

For any $(\G'\sub\G'')$ in $\bar x_\G$
and any $(\G'_1\sub\G''_1)$ in $\bX_{\G''/\G'}$ 
we define $\ti\G'_1,\ti\G''_1$ to be the inverse images of
$\G'_1,\G''_1$ under the quotient map $\G''@>>>\G''/\G'$.

The pairs $(\ti\G'_1\sub\ti\G''_1)$ thus associated to various
$(\G'\sub\G'')$ in $\bar x_\G$ 
and any $(\G'_1\sub\G''_1)$ in $\bX_{\G''/\G'}$ 
form a subset $(\bX_\G)_0$ of $\cz_\G$.
If $\G$ is not anomalous and not terminal we define 
$\bX_\G=(\bX_G)_0\sqc\{(S_1\sub\G)\}$.
If $\G$ is anomalous or terminal
we set $\bX_\G=X_\G\sqc\{(S_1\sub S_1)\}$.

This completes the inductive definition of $\bX_\G$.

Note that for any $(\G'_1\sub\G''_1)$ in $\bX_\G$,
$\G''_1/\G'_1$ is isomorphic to a product of groups in $\AA$.

\subhead 3.3\endsubhead
Combining \cite{L22}, \cite{L23} with the results of this paper (but with the
sets $x_{\G_c},X_{\G_c}$ enlarged to sets $\bar x_{\G_c},\bX_{\G_c}$
one can obtain an indexing of the set of unipotent character sheaves
of $G$ corresponding to a family $c$ extending the indexing of the
representations of $W$ in $c$.

\head 4. Precuspidal families\endhead
\subhead 4.1\endsubhead
Let $\{s_i;i\in I\}$ be the set of simple reflections of $W$.
For any $I'\sub I$ let $W_{I'}$ be the subgroup of $W$ generated by
$\{s_i;i\in I'\}$; this is again a Weyl group whose set of simple
reflections is $\{s_i;i\in I'\}$.

We say that an irreducible Weyl group $W'$ is {\it non-terminal}
(resp. {\it terminal}) if we can (resp. cannot) find an irreducible
$W,I$ as above and $I'\sneq I$ such that $W'=W_{I'}$; the condition
that $W'$ is terminal is equivalent to the condition that $W'$ is of
type $G_2,F_4$ or $E_8$. 

\subhead 4.2\endsubhead
Let $\sgn=\sgn_W\in\hW$ be the sign representation of $W$.

If $W$ is a
product $W_1\T W_2\T\do\T W_k$ where $W_j$ are irreducible
Weyl groups, then we have a bijection
$\Ph(W_1)\T\Ph(W_2)\T\do\T\Ph(W_k)@>\si>>\Ph(W)$
given by
$$(c_1,c_2,\do,c_k)\m\{E\in\hW;E=E_1\bxt E_2\bxt\do\bxt E_k
\text{ with }E_1\in c_1,\do,E_k\in c_k\}.\tag a$$

When $(c_1,c_2,\do,c_k)\m c$ as in (a), we say that
$c_1,c_2,\do,c_k$ are the {\it components} of $c$.

In \cite{L84} to any $c\in\Ph(W)$ we have attached a
an imbedding
$\io_c:c@>>>M(\G_c)$ (notation of 2.4) with image $M(\G_c)_0$.

Let $E\m a_E$ be the function $\hW@>>>\NN$
defined in \cite{L84, 4.1}. It is
known \cite{L84} that $E\m a_E$ is constant on each $c\in\Ph(W)$.

If $I'\sub I,E\in\hW,E'\in\hW_{I'}$, we denote by $<E',E>_W$
the multiplicity of $E'$ in the restriction of $E$ to
$W_{I'}$. For $E'\in\hW_{I'}$ we set
$$J_{W_{I'}}^W(E')=\sum_{E\in\hW;a_{E'}=a_E}<E',E>_WE\in\car(W).$$
This extends to a linear map 
$J_{W_{I'}}^W:\car(W_{I'})@>>>\car(W)$.
This restricts for any $c'\in\Ph(W_{I'})$ to a linear map
$\car_{c'}@>>>\car_c$
where $c$ is a well defined family of $W$ denoted by
$J_{W_{I'}}^W(c')$.

\subhead 4.3\endsubhead
Assuming that $W$ is irreducible we describe the group $\G_c$
attached to $c\in\Ph(W)$. It is an object of $\AA$.

If $W$ is of type $A_n,n\ge1$ we have $|c|=1,\G_c=S_1$.

If $W$ is of type $B_n$ or $C_n$, $n\ge2$, we have $\G_c=V_D^1$ for
some even $D\ge0$.

If $W$ is of type $D_n,n\ge4$ we have $\G_c=V'_D{}^1$ for
some odd $D\ge0$.

If $W$ is of exceptional type, then we are in one of the
following cases.

$|c|=1,\G_c=S_1$;

$|c|=2,\G_c=S'_2$;

$|c|=3,\G_c=S_2$;

$|c|=4,\G_c=S'_3$;

$|c|=5,\G_c=S_3$;

$|c|=11,\G_c=S_4$;

$|c|=17,\G_c=S_5$.

\subhead 4.4\endsubhead
Note that:

(a) {\it if $c\in\Ph(W)$, then
$c\ot\sgn_W:=\{E\ot\sgn_W;E\in c\}\in\Ph(W)$.}
\nl
We have $\G_{c\ot\sgn_W}=\G_c$.

Let $c\in\Ph(W)$. Following \cite{L82} we say that $c$ is
{\it smoothly induced} if there exist $I'\sneq I$ and
$c'\in\Ph(W_{I'})$ such that $E'\m J_{W_{I'}}^W(E')$ (notation
of 4.2) is a bijection $c'@>\si>>c$; in this case we have
$\G_c=\G_{c'}$.

We say that $c$ is {\it cuspidal} if $c$ is not smoothly induced and
$c\ot\sgn_W$ (see (a)) is not smoothly induced.
If $c\in\Ph(W)$ is cuspidal we set $\g(c)=0$ if $|I|$ is even,
$\g(c)=1$ if $|I|$ is odd.

The following result is implicit in \cite{L84}.

(b) $\Ph(W)$ {\it contains at most one cuspidal family.}
\nl
When $W$ is irreducible and $c\in\Ph(W)$, we say that $c$ is
{\it terminal} if $c$ is cuspidal and $W$ is terminal that is of type
$G_2,F_4$ or $E_8$.

When $W$ is irreducible and $c\in\Ph(W)$, we say that $c$ is
{\it anomalous} if $|c|=2$.

\subhead 4.5\endsubhead
Assume that $W$ is irreducible and $c\in\Ph(W)$ is cuspidal
and not anomalous.
Let $\Si'_c$ be the set of all pairs $(I',c')$ where $I'\sub I$,
$|I'|=|I|-1$, $c'\in\Ph(W_{I'})$ is not smoothly induced and
$c=J_{W_{I'}}^W(c')$.       

Let $\Si_c$ be the set of all $(I',c')\in\Si'_c$ with the following
property: if $c''$ is a component of $c'$ (see 4.2), then either
$c''$ is non-cuspidal or $c''$ is cuspidal with $\g(c'')=\g(c)$. 

We say that the families $c'$ which appear in some pair
$(I',c')\in\Si_c$ are the {\it precuspidal} families associated to $c$.

We describe explicitly in each case the sets $\Si_c,\Si'_c$ attached
to $c$. We will specify a family $c'$ by specifying the
corresponding special representation $E_{c'}$ in $c'$ of the Weyl group.
(The notation is that of \cite{L84a}; in particular for classical types
we use the symbol notation.) We will also specify the group $\G_c$.

(a) If $W$ is of type $B_2=C_2$, we have $E_c=[0.2.4//.1.3.]$  
$\G_c=V_2^1$ and
$\Si_c=\Si'_c$ consists of the two pairs $(I',\sgn)$ with $|I'|=1$.

(b) If $W$ is of type $B_{k^2+k}=C_{k^2+k}$ with $k\ge2$, then
$$E_c=[0.2.4.\do.(2k)//.1.3.5.\do.(2k-1).],$$ $\G_c=V_{2k}$
and $\Si_c=\Si'_c$ consists
of $2k$ pairs $(I',c')$ where $W_{I'}$  are of type
$$B_{k^2+k-1},B_{k^2+k-2}A_1,B_{k^2+k-3}A_2,\do,B_{k^2-k}A_{2k-1};$$
the corresponding $E_{c'}$ are
$$[0.2.4.\do.(2k-2).(2k-1)//.1.3.5.\do.(2k-1).],$$
$$[0.2.4.\do.(2k-2).(2k-1)//.1.3.5.\do.(2k-3).(2k-2).]\bxt\sgn,$$
$$[0.2.4.\do.(2k-4).(2k-3).(2k-1)//.1.3.5.\do.(2k-3).(2k-2).]\bxt\sgn,$$
$$\do$$
$$[0.1.3.5.\do.(2k-1)//.1.2.4.\do.(2k-2).]\bxt\sgn,$$
$$[0.1.3.5.\do.(2k-1)//.0.2.4.\do.(2k-2).]\bxt\sgn.$$
(Note that for $(I',c')$ such that $W_{I'}$ is of type 
$B_{k^2-k}A_{2k-1}$, the $B_{k^2-k}$-component $c''$ of $c'$ is
cuspidal with $\g(c'')=\g(c)=0$.)

(c) If $W$ is of type $D_4$, then $E_c=[0.2.//.1.3]$, $\G_c=V'_3{}^1$
and $\Si_c=\Si'_c$ consists
of the three $(I',\sgn)$ with $W_{I'}$ of type $A_2$ and of the
unique $(I',\sgn)$ with $W_{I'}$ of type $A_1A_1A_1$.

(d) If $W$ is of type $D_{k^2}$ with $k\ge3$,
then $$E_c=[0.2.4.\do.(2k).//.1.3.5.\do.(2k+1)],$$
$\G_c=V'_{2k-1}{}^1$
and $\Si_c$ consists of $2k-2$ pairs $(I',c')$ where $W_{I'}$ are of
type
$$D_{k^2-1},D_{k^2-2}A_1,D_{k^2-3}A_2,\do,D_{k^2-2k+2}A_{2k-3};$$
the corresponding $E_{c'}$ are
$$[0.2.4.\do.(2k-2).(2k)//.1.3.5.\do.(2k-1).(2k)]\bxt\sgn,$$
$$[0.2.4.\do.(2k-2).(2k-1)//.1.3.5.\do.(2k-1).(2k)]\bxt\sgn,$$
$$[0.2.4.\do.(2k-2).(2k-1)//.1.3.5.\do.(2k-3).(2k-2).(2k)]\bxt\sgn,$$
$$\do$$
$$[0.1.3.5.\do.(2k-1)//.1.2.4.\do.(2k)]\bxt\sgn.$$
$\Si'_c$ is the union of $\Si_c$ and one other pair $(I',c')$ with
$W_{I'}$ of type $D_{(k-1)^2}A_{2k-2}$ and $E_{c'}$ given by
$$[0.1.3.5.\do.(2k-1)//.0.2.4.\do.(2k)]\bxt\sgn.$$
(Note that for the last $(I',c')$, the $D_{(k-1)^2}$-component $c''$ of $c'$ is cuspidal with $\g(c'')\ne\g(c)$.)

(e) If $W$ is of type $E_6$, then $E_c=80_s$, $\G_c=S_3$ and
$\Si_c=\Si'_c$ consists of the two $(I',c')$ with $W_{I'}$ of type $D_5$,
$E_{c'}=[.1.2.4//0.1.3.]$ and of the unique $(I',\sgn)$ with $W_{I'}$
of type $A_2A_2A_1$.

(f) If $W$ is of type $E_8$, then $E_c=4480_y$, $\G_c=S_5$ and
$\Si_c=\Si'_c$ consists of the five $(I',c')$ with $W_{I'}$ of type
$E_7,E_6A_1,D_7,D_5A_2,A_4A_3$ and $E_{c'}$ given respectively by
$$315_a,30'_p\bxt\sgn,[.1.2.3.5//0.1.3.4.],[.1.2.3.4//0.1.2.4.]\ot\sgn,
\sgn.$$
 
(g) If $W$ is of type $F_4$, then $E_c=12_1$, $\G_c=S_4$ and
$\Si_c=\Si'_c$ consists of the two $(I',c')$ with $W_{I'}$ of type
$B_3,C_3$, $E_{c'}=[0.1.3//.1.2.]$, and the two $(I',\sgn)$ with
$W_{I'}$ of type $A_2\T A_1$.

(h) If $W$ is of type $G_2$, then $E_c=V$, $\G_c=S'_3$ and
$\Si_c=\Si'_c$ consists of the two $(I',\sgn)$ such that $|I'|=1$.

\subhead 4.6\endsubhead
Assume that $W$ is irreducible and that $c\in\Ph(W)$ is cuspidal
and anomalous. (Thus $W$ is of type $E_7$.) We have $E_c=512'_a$,
$\G_c=S'_2$. Let $\Si_c=\Si'_c$ be the set of all pairs $(I',c')$ where
$I'\sub I$ is such that $W_{I'}$ is of type $A_4\T A_1$ and
$c'\in\Ph(W_{I'})$ is given by $\sgn$. Note that $c'$ is not smoothly
induced and
$c=J_{W_{I'}}^W(c')$. (But unlike in 4.5 we do not have $|I'|=|I|-1$.)   We again say that $c'$ is a precuspidal family attached to $c$.

\subhead 4.7\endsubhead
For any finite group $\G$ we set $\un\G=\sum_x(x,1)\in\CC[M(\G)]$ where
the sum is taken over a set of representatives for the conjugacy classes
in $\G$. For any subgroup $\G''$ of $\G$ we set
$$[\G'']_\G=\sum_{(x,\s)\in M(\G)}m_{x,\s}(x,\s)\in\CC[M(\G)]$$
where $m_{x,\s}$ is the multiplicity of $\s$ in the permutation
$Z_\G(x)$-module defined by the fixed point set 
$\{g\G'';xg\G''=g\G''\}$ of $x$ on $\G/\G''$. Note that
$[\G]_\G=\un\G$,
$$[\{1\}]_\G=\sum_{(1,\s)\in M(\G)}\dim(\s)(1,\s)\in\CC[M(\G)].$$
From the definitions, for $(\G'\sub\G'')\in\cz_\G$ we have
$$\ss_{\G',\G''}(\un{\G''/\G'})=[\G'']_\G.\tag a$$
(Notation of 2.4).

\subhead 4.8\endsubhead
We assume that $W$ is irreducible and that $c\in\Ph(W)$ is cuspidal. We
fix $(I',c')\in\Si_c$. The linear map $J_{W_{I'}}^W:\car_{c'}@>>>\car_c$
(see 4.2) can be viewed as a linear map
$\CC[M(\G_{c'})_0]@>>>\CC[M(\G_c)_0]$ via the bijection
$$c@>\si>>M(\G_c)_0\tag a$$
induced by $\io_c$ in 4.2 and the analogous
bijection $c'@>\si>>M(\G_{c'})_0$.
(Note that (a) gives rise to an identification

(b) $\car_c=\CC[M(\G_c)_0]$.
\nl
With notation in 4.7 we have:
$$J_{W_{I'}}^W(\un{\G_{c'}})=[\G'']_{\G_c}\tag c$$
for a subgroup $\G''$ of $\G_c$ (well defined up to conjugacy) which
is endowed with a surjective homomorphism $\G''@>>>\G_{c'}$.
This follows from the explicit description of the $J$-induction given in
\cite{L84, (4.5.4), (4.6.5), 4.10, 4.11, 4.12, 4.13}.
(See also \S5 for an alternative approach.)
Let $\G'$ be the kernel of $\G''@>>>\G_{c'}$.
The following result follows from the examination of the various cases.

\proclaim{Theorem 4.9} We assume that $W$ is irreducible and that
$c\in\Ph(W)$ is cuspidal. The pairs $(\G'\sub\G'')\in\cz_{\G_c}$
associated in 4.8 to the various $(I',c')\in\Si_c$ form precisely the
set $x_{\G_c}$ (see 2.2).
\endproclaim
More precisely, $(I'_1,c'_1),(I'_2,c'_2)$ in $\Si_c$ give rise to the
same $(\G'\sub\G'')$ if and only if conjugation by some
element of $W$ carries $(I'_1,c'_1)$ to $(I'_2,c'_2)$.

\head 5. Relation to unipotent classes\endhead
\subhead 5.1\endsubhead
Let $P_\emp$ be a fixed Borel subgroup of $G$ in 0.1. For any $I'\sub I$
let $P_{I'}$ be the parabolic subgroup of type $I'$ containing $P_\emp$;
let $L_{I'}$ be the reductive quotient of $P_{I'}$ and let
$p_{I'}:P_{I'}@>>>L_{I'}$ be the obvious surjective map. If
$c\in\Ph(W)$, then $c$ contains a unique special representation, see
\cite{L84}, and that representation is associated to a special 
unipotent class $\bx{c}$ of $G$ as in \cite{L84, (13.1.1)}. For
$u\in\bx{c}$ let $\bA_G(u)$ be the quotient of $A_G(u)$ (see 0.2)
defined in \cite{L84, 13.1}; it is known that $\bA_G(u)$ is isomorphic
to $\G_c$.

\subhead 5.2\endsubhead
We now assume that $G$ is simple and $c\in\Ph(W)$ is cuspidal.
For any $(I',c')\in\Si_c$, the special unipotent class $\bx{c}$ of $G$ is
induced in the sense of \cite{LS79} from the special unipotent class
$\bx{c'}$ of $L_{I'}$; thus $\bx{c}\cap p_{I'}\i\bx{c'}$ is dense in
$p_{I'}\i\bx{c'}$. Let $u\in\bx{c}\cap p_{I'}\i\bx{c'}$ and let
$u'=p_{I'}(u)$. Let $\G''=A_{P_{I'}}(u)$ (see 0.2). This is a
subgroup of $A_G(u)=\bA_G(u)=\G_c$, see \cite{LS79, 1.3(d)}. Let $\G'$
be the kernel of the homomorphism $A_{P_{I'}}(u)@>>>A_{L_{I'}}(u')$
induced by $p_{I'}$; this homomorphism is surjective by
\cite{LS79, 1.5}. Note that $\G'$ is a normal subgroup of $\G''$. One
can verify that

(a) the pair $(\G'\sub\G'')$ just defined is up to $\G_c$-conjugacy
the same as the pair $(\G'\sub\G'')$ associated to $(I',c')$ in 4.7.
\nl
We thus obtain an an alternative description (in terms of
unipotent elements) of the set $x_{\G_c}$.

\subhead 5.3\endsubhead
Let $V$ be a $\CC$-vector space of dimension $24$ with a nondegenerate
symplectic form $(,)$. In this subsection we assume that $G$ is the
symplectic group of $V,(,)$ modulo its centre. We fix a basis
$$\align&\{e^1_1,e^1_2,f^2_1,f^2_2\}\sqc\\&
\{e^2_1,e^2_2,e^2_3,e^2_4,f^2_1,f^2_2,f^2_3,f^2_4\}\sqc\\&
\{e^3_1,e^3_2,e^3_3,e^3_4,e^3_5,e^3_6,\\&
f^3_1,f^3_2,f^3_3,f^3_4,f^3_5,f^3_6\}\endalign$$
of $V$ such that any two basis elements have $(,)=0$ except for
$(e^j_a,f^j_b)=(-1)^a\d_{a+b,2j+1}$
$(f^j_b,e^j_a)=-(-1)^a\d_{a+b,2j+1}$
for $j=1,2,3$.
Let $N:V@>>>V$ be the nilpotent linear map such that
$$e^1_1\m e^1_2\m0,f^2_1\m f^2_2\m0,
e^2_1\m e^2_2\m e^2_3\m e^2_4\m0,f^2_1\m f^2_2\m f^2_3\m f^2_4\m0,$$
$$e^3_1\m e^3_2\m e^3_3\m e^3_4\m e^3_5\m e^3_6\m0,
f^3_1\m f^3_2\m f^3_3\m f^3_4\m f^3_5\m f^3_6\m0.$$

Note that $N$ has Jordan blocks of sizes $6,6,4,4,2,2$ and that
$(Nv,v')+(v,Nv')=0$ for any $v,v'$ in $V$.
Consider the following subspaces of $V$:

$V_6=\text{span}\{e^3_6,f^3_6,e^2_4,f^2_4,e^1_2,f^1_2\}$,

$V_5=\text{span}\{e^3_6,f^3_6,e^2_4,f^2_4,e^1_2\}$,

$V_4=\text{span}\{e^3_6,f^3_6,e^2_4,f^2_4\}$,

$V_3=\text{span}\{e^3_6,f^3_6,e^2_4\}$,

$V_2=\text{span}\{e^3_6,f^3_6\}$,

$V_1=\text{span}\{e^3_6\}$.

Let $V_h^\pe:=\{v\in V;(x,V_h)=0\}$.
Note that for $h\in[1,6]$, we have $(V_h,V_h)=0$, $N|_{V_h}=0$ and $N$
induces on $V_h^\pe/V_h$ a nilpotent linear map with Jordan blocks of
sizes

$4,4,2,2$   if $h=6$,

$4,4,2,2,1,1$  if $h=5$,

$4,4,2,2,2,2$   if $h=4$,

$4,4,3,3,2,2$  if $h=3$,

$4,4,4,4,2,2$   if $h=2$,

$5,5,4,4,2,2$   if $h=1$.

Now $u=\exp(N)$ can be viewed as a unipotent element of $G$ as in 5.2
and for any $h\in [1,6]$, $u'=(1|V_h)\T\exp(N|V_h^\pe/V_h)$ can be
viewed as a unipotent element of a group $L_{I'}$ in 5.1;
moreover, $u'$ is as in 5.2.  
In particular, the conjugacy class of $u$ in $G$ is induced in the sense
of \cite{LS79} from the conjugacy class of $u'$ in $L$.

For $j=1,2,3$ we define $T_j:V@>>>V$ by

$T_j(e^j_a)=f^j_a$, $T_j(f^j_a)=e^j_a$,

$T_j(e^i_a)=e^i_a$, $T_j(f^i_a)=f^i_a$ for $i\ne j$.
\nl
Then $A_G(u)=\bA_G(u)$ can be identified with the commutative $2$-group
$\ct$ consisting of $T_1^{c_1}T_2^{c_2}T_3^{c_3}$
where $c_j\in\{0,1\}$ for $j=1,2,3$. For $h\in[1,6]$ we set
$\ct_h=\{\t\in\ct;\t(V_h)=V_h\}$.

When $h$ is even we have $\ct_h=\ct$; when $h$ is odd, $\ct_h$
consists of all $T_1^{c_1}T_2^{c_2}T_3^{c_3}\in\ct$ such that $c_3=0$
(if $h=1$), $c_2=0$ (if $h=3$), $c_1=0$ (if $h=5$). It follows that
5.2(a) holds in this case. 

\head 6. The basis $\fF_c$ of $\car c$\endhead
\subhead 6.1\endsubhead
In this subsection we assume that $W=S_n$ where $n\ge1$ (see 2.1).
Let $Part(n)$ be the set of partitions of $n$ that it symbols
$a_*=a_1a_2a_3\do$ with $a_i\in\NN$ ($a_i=0$ for large $i$) such that
$a_1\ge a_2\ge a_3\ge\do$ and $a_1+a_2+a_3+\do=n$.
For $a_*,a'_*$ in $Part(n)$ we write $a_*\ge a'_*$ whenever
$a_1\ge a'_1$, $a_1+a_2\ge a'_1+a'_2$, etc; this is a partial order
on $Part(n)$. For example if $n=3$ we have $3\ge21\ge111$; if $n=4$ we
have $4\ge31\ge22\he211\ge1111$; if $n=5$ we have
$5\ge41\ge32\ge311\ge221\ge2111\ge11111$
(we omit writing zeros in $a_*$) so that in these cases we actually
have a total order. (This is not so for $n\ge6$.)

For any $a_*\in Part(n)$ we denote by $S(a_*)$ the subgroup of  $S_n$
consisting of all permutations of $[1,n]$ which preserve each of the
subsets

$\{1,2,\do,a_1\},\{1,2,\do,a_1+a_2\},\{1,2,\do,a_1+a_2+a_3\}$,
\nl
etc. Let $r_{a_*}=\Ind_{S(a_*)}^{S_n}(1)$. We can regard $r_{a_*}$
as an element of $\car(W)=\car(S_n)$ (see 4.2). According to Frobenius,
there is a unique way to index the elements of $\hW$ (see 0.1) as
$\{E_{a_*}; a_*\in Part(n)\}$ such that for any $a_*$ we have
$r_{a_*}-E_{a_*}\in\sum_{a'_*>a_*}\NN E_{a'_*}$.

For $n=3$ we have
$r_3=E_3, r_{21}=E_{21}+E_3,r_{111}=E_{111}+2E_{21}+E_3.$

For $n=4$ we have

$r_4=E_4, r_{31}=E_{31}+E_4,r_{22}=E_{22}+E_{31}+E_4,$

$r_{211}=E_{211}+E_{22}+2E_{31}+E_4,$

$r_{1111}=E_{1111}+3E_{211}+2E_{22}+3E_{31}+E_4.$

For $n=5$ we have

$r_5=E_5, r_{41}=E_{41}+E_5,r_{32}=E_{32}+E_{41}+E_5,$

$r_{311}=E_{311}+E_{32}+2E_{41}+E_5,
r_{221}=E_{221}+E_{311}+2E_{32}+2E_{41}+E_5,$

$r_{2111}=E_{2111}+2E_{221}+3E_{311}+3E_{32}+3E_{41}+E_5$,

$r_{11111}=E_{11111}+4E_{2111}+5E_{221}+6E_{311}+5E_{32}+4E_{41}+E_5.$

\subhead 6.2\endsubhead
We return to a general $W$. In this subsection we assume that
$W$ is irreducible and $c\in\Ph(W)$ is terminal (see 4.4) hence
cuspidal. Then $\G_c=S_n$ where $n\in\{3,4,5\}$. (Here we identify
$S'_3=S_3$ as groups.)
If $E$ is an irreducible representation of $\G_c$, the element
$(1,E)\in M(\G_c)$ belongs to $M(\G_c)_0$ (see 4.2) if and only if
$E\ne E^0$ where $E^0$ is the sign representation of $S_n$. We define
a linear map $\t:\car(S_n)@>>>\CC[M(\G_c)_0]=\car_c$ by
$E_{a_*}\m(1,E_{a_*})$ if $E_{a_*}\ne E^0$, $E_{a_*}\m0$ if
$E_{a_*}=E^0$. (Notation of 6.1).
Let $\car_c^!$ be the subset of $\car_c$ consisting of the elements
$\t(r_{a_*})$ (notation of 6.1) for various $a_*\in Part(n)$ such that
$E_{a_*}\ne E^0$ (that is $a_*$ gas some part $\ge2$). This is a
linearly independent subset.

\subhead 6.3\endsubhead
For any $c\in\Ph(W)$ we define a subset $\fF_c$ of $\car_c$ by
induction on $|I|$. If $|I|=0$, $\fF_c$ consists of
$1\in\car_c=\CC$. Assume now that $|I|>0$.

If there exists $I'\sneq I$ and $c'\in\Ph(W')$ such that $c$ is
smoothly induced (see 4.4) from $I',c'$ then $J_{W_{I'}}^W$ defines
an isomorphism
$\car_{c'}@>>>\car_c$; we define $\fF_c$ to be the set of elements of
$\car_c$ obtained by applying this isomorphism to the elements
in $\fF_{c'}$. One shows that that

(a) {\it this is independent of choices.}
\nl
If there exists $I'\sneq I$ and $c'\in\Ph(W')$ such that
$c\ot\sgn$ is smoothly induced from $I',c'$ then
$\fF_{c\ot\sgn_W}$ is defined by the previous paragraph.
Under the obvious isomorphism $\car_c@>>>\car_{c\ot\sgn_W}$,  
$\fF_{c\ot\sgn_W}$ becomes a subset $\fF_c$ of $\car_c$.

If $c$ is not as in the previous two paragraphs, then $c$ is cuspidal.
If $W$ is a product $W_1\T W_2\T\do\T W_k$ where $W_j$ are irreducible
Weyl groups with $k\ge2$
and $(c_1,c_2,\do,c_k)\m c$ are as in 4.2(a) (and are necessarily
cuspidal) then
we can identify $\car_{c_1}\ot\car_{c_2}\ot\do\ot\car_{c_k}=\car_c$;
we define $\fF_c$ to be the set of elements
$\cx_1\ot\cx_2\ot\do\ot\cx_k$ where $\cx_j\in\fF_{c_j}$ for all $j$.

Thus we can assume that $c$ is cuspidal and $W$ is irreducible,
so that the set $\Si_c$ is defined as in 4.5, 4.6. We say that
an element $X\in\car_c$ is in $\fF_c$ if one of (i),(ii),(iii) below
holds:

(i) $c$ is not terminal and $X$ is the special representation in $c$;

(ii) $c$ is terminal and $X\in\car_c^!$ (see 6.2);

(iii) there exists $(I',c')\in\Si_c$ and $X'\in\fF_{c'}$ such that $X$
is the image of $X'$ under the linear map
$J_{W_{I'}}^W:\car(W_{I'})@>>>\car(W)$.
\nl
This completes the inductive definition of $\fF_c$.

\subhead 6.4\endsubhead
In the following theorem we assume that $W$ is irreducible and that
$c\in\Ph(W)$. Via the identification 4.8(b), the basis of
$\CC[M(\G_c)_0]$ described in 2.5 becomes a basis of $\car_c$.

\proclaim{Theorem 6.5} This basis of $\car_c$ coincides with
$\fF_c$. In particular, $\fF_c$ is a basis of $\car_c$.
\endproclaim
This follows again from the explicit description of the $J$-induction
in \cite{L84} referred to in 4.8. We use the fact that (assuming that
$c$ is cuspidal) the maps $\ss_{\G',\G''}$ (see 2.4) are very closely
connected to $J$-induction from parabolic subgroups of $W$. (This
connection has already been pointed out in type $E_8$
in \cite{L84, p.311}.)

\enddocument